\documentclass[11pt,final]{amsart}
\usepackage{amsaddr}
\usepackage{mathpazo}
\usepackage{parskip}
\usepackage[pagewise]{lineno}

\usepackage[centering]{geometry}
\usepackage{enumerate}
\usepackage{amsmath}
\usepackage{amsthm}
\usepackage{amssymb}
\usepackage{amscd}
\usepackage{amsfonts}
\usepackage{amsbsy}
\usepackage{graphicx}
\usepackage{caption}
\usepackage{subcaption}
\usepackage[dvips]{psfrag}
\usepackage{array}
\usepackage{epsfig}
\usepackage{url}
\usepackage{overpic}
\usepackage{epstopdf}
\usepackage[normalem]{ulem}
\usepackage{tikz}
\usepackage{overpic}
\usepackage{multirow}
\usetikzlibrary{arrows}
\usepackage{overpic}
\usepackage{float}
\usepackage[pdftex]{hyperref}

\usepackage[utf8]{inputenc} 
\usepackage[T1]{fontenc}    

\usepackage{url}            
\usepackage{booktabs}       
\usepackage{amsfonts}       
\usepackage{nicefrac}       
\usepackage{microtype}      
\usepackage{lipsum}
\usepackage{ dsfont }
\usepackage[utf8]{inputenc}
\usepackage{caption}
\usepackage{amsmath}
\usepackage{amssymb}
\usepackage{pgfpages}
\usepackage{amsthm}
\usepackage{array}
\usepackage{multirow}
\usepackage{mathrsfs}

\usepackage{ragged2e}
\usepackage{graphicx}
\usepackage{graphics}

\newcommand{\C}{\ensuremath{\mathbb{C}}}

\usepackage[toc,page]{appendix}
\usepackage{tikz-cd}
\theoremstyle{definition}
\theoremstyle{plain}
\newtheorem{theorem}{Theorem}
\newtheorem{remark}{Remark}

\newtheorem{lema}{Lemma}
\newtheorem{exemplo}{Example}

\newtheorem*{bezout}{Theorem}
\theoremstyle{remark}

\usepackage{xargs}

\usepackage{todonotes}
\usepackage{pdfcomment}

\begin{document}

\title[Limit cycles of PSVF on torus]{Limit cycles of piecewise smooth vector fields on torus with non-regular switching manifold}

\author[T. S. Oliveira, R. M. Martins]{Thaylon Souza de Oliveira$^{\ast}$ and Ricardo Miranda Martins}

\address{Universidade Estadual de Campinas (UNICAMP)\\ Department of Mathematics, Institute of Mathematics, Statistics and Scientific Computing. Campinas/SP, Brazil, 13083-859}
\thanks{$^\ast$Corresponding author.}

\email{thaylon@ime.unicamp.br, rmiranda@unicamp.br}



\begin{abstract}
We study the maximum number of crossing limit cycles in piecewise‐smooth vector fields on the two‐dimensional torus, where the discontinuity set is the boundary of the fundamental square. Under the assumption of a polynomial first integral of degree \(n\), we apply algebraic curve‐intersection methods to obtain sharp upper bounds on cycle counts. In the quadratic case (\(n=2\)), we prove at most one cycle of type \emph{\textbf{aa}} or \(\emph{\textbf{bb}}\) and at most two of type \(\emph{\textbf{aba}}\) or \(\emph{\textbf{bab}}\), and we give explicit parameter criteria for their existence. For general \(n\), we show that cycles of type \(\emph{\textbf{bb}}\) are bounded by \(n-1\) and those of type \(\emph{\textbf{aba}}\) by \(n(n-1)\). Finally, we construct concrete examples achieving these bounds, demonstrating their optimality.
\end{abstract}

\maketitle


\section{Introduction and statement of the main result}

The theory of piecewise smooth vector fields (PSVF) has attracted growing attention in recent years due to its close connection with various fields of applied sciences, such as mathematics, physics, and engineering. Recent studies discuss the importance of PSVF, including models that arise in control theory (see, for instance, \cite{makarenkov2012dynamics,teixeira2022perturbation,bernardo2008piecewise, simpson2010bifurcations}). In general, PSVF consist of smooth differential systems defined in distinct open regions of a manifold, where the boundary between these regions, which separates the different smooth vector fields, is known as the switching manifold (or discontinuity manifold).

Research in PSVF has made significant progress in characterizing local behaviors near the switching manifold, particularly in defining singularities and trajectories in planar cases \cite{guardia2011generic}. Furthermore, topics such as bifurcations and stability have been explored in both general and symmetric systems \cite{broucke2001structural,jacquemard2011piecewise,jacquemard2013stability,kuznetsov2003one,makarenkov2012dynamics,teixeira1990stability,teixeira1993generic}. These theoretical advancements find practical applications in control and relay systems, where specific phenomena such as sliding motion are observed. Sliding motion occurs when trajectories reach the switching manifold and travel along it until reaching a boundary, forming what is known as the sliding region. This behavior, characteristic of relay control and dry friction systems \cite{bernardo2008piecewise,di2011two,chillingworth2002discontinuity,jacquemard2011piecewise,jacquemard2013stability,jacquemard2012coupled,jeffrey2009two,makarenkov2012dynamics}, can manifest in three main types of regions on the discontinuity manifold: crossing, sliding and escape, which together with the fold regions, describe the range of possible dynamics along the manifold.

It is natural to study limit cycles in piecewise smooth vector fields, due to its theoretical importance and also for applications \cite{filippov1960differential, bernardo2008piecewise,di2011two, simpson2010bifurcations}. By understanding these limit cycles, we gain insight into the stability and periodicity within discontinuous dynamical systems \cite{makarenkov2012dynamics, teixeira1993generic}. 
%
For instance, in \cite{llibre2018limit} the dynamics of PSVF on a two-dimensional torus is explored, and the authors prove that some linear and Riccati-like vector fields generates limit cycles, with the potential to produce an arbitrarily large number of limit cycles. Similarly, in \cite{martins2019chaotic}, authors  investigated chaotic behavior in PSVF on the torus and sphere, demonstrating how discontinuities create non-deterministic, dense trajectories under certain bifurcation conditions. These findings emphasize the relevance of studying PSVF on the torus to understand systems with cyclic and chaotic responses to switching boundaries.



We focus on PSVF defined on the two‑dimensional torus \(\mathbb{T}^2\).  Concretely, we model \(\mathbb{T}^2\) as the quotient of the unit square \([0,1]^2\) with opposite edges identified.  The discontinuity set is precisely the square’s boundary: trajectories crossing an edge reappear on the ``opposite side''.

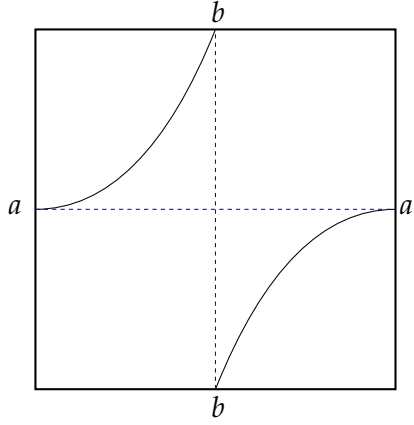
\begin{figure}[h]
  \centering
  \begin{overpic}[height=5cm,unit=1mm]{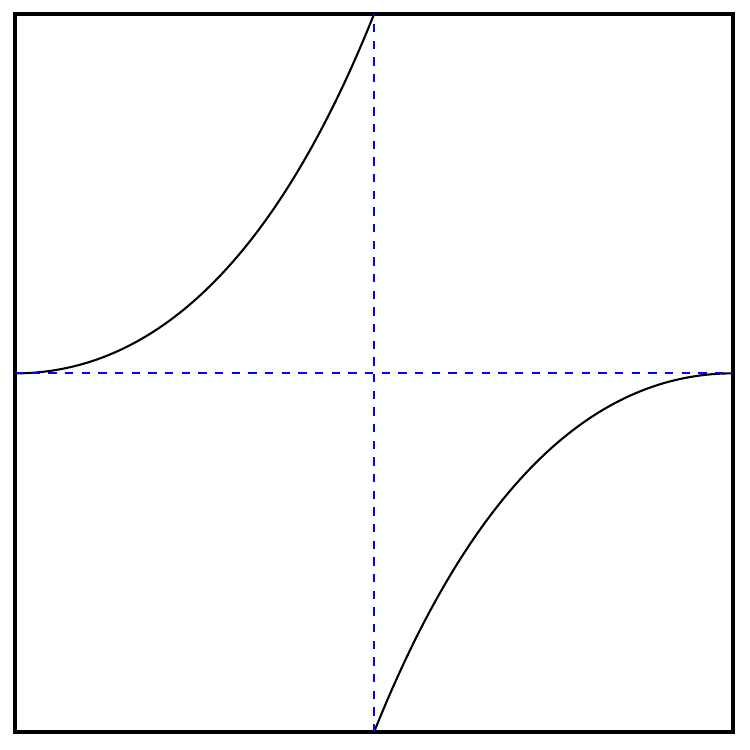}
    \put(-4,49){\(\!a\)}
    \put(100,49){\(\!a\)}
    \put(50,-5){\(\!b\)}
    \put(50,100){\(\!b\)}
  \end{overpic}
  \caption{Symbolic representation of crossing cycles: \(\emph{\textbf{aba}}\) on the torus.}
  \label{fig:symbolic}
\end{figure}

Our main interest lies in \emph{crossing limit cycles}—closed orbits that intersect the switching manifold transversally, and are unique with this property in a neighborhood. We encode such cycles by symbolic sequences of horizontal (\textbf{a}) and vertical (\textbf{b}) loops: for instance, a cycle of type \emph{\textbf{aa}} corresponds to two consecutive horizontal loops, while a cycle of type \emph{\textbf{bb}} corresponds to two consecutive vertical loops. Similarly, cycles of type \emph{\textbf{aba}} and \emph{\textbf{bab}} describe a horizontal loop, a vertical loop, and the return horizontal loop, and the reverse sequence, respectively. Figure~\ref{fig:symbolic} illustrates these conventions.

In this work, we assume the existence of a polynomial first integral \(H(x,y)\) of degree \(n\) for our systems. Such integrals are fundamental tools in the analysis of dynamical systems, as they constrain trajectories to lie on level curves of \(H\). In our setting, the closing conditions across the switching manifold translate into algebraic relations, which we refer to as ``closing equations.''

By applying algebraic curve–intersection techniques (notably projective Bézout’s theorem) and verifying transversality hypotheses, we derive sharp upper bounds on the number of crossing cycles of types \(\emph{\textbf{aa}}\), \(\emph{\textbf{bb}}\), \(\emph{\textbf{aba}}\), and \(\emph{\textbf{bab}}\). In particular, we prove that:
\begin{itemize}
  \item for quadratic first integrals (\(n=2\)), at most one cycle of type \emph{\textbf{aa}} or \(\emph{\textbf{bb}}\) and at most two of type \(\emph{\textbf{aba}}\)/\(\emph{\textbf{bab}}\) can occur.
  \item for general \(n\), the number of \(\emph{\textbf{bb}}\)-type cycles is bounded by \(n-1\), and the number of \(\emph{\textbf{aba}}\)-type cycles by \(n(n-1)\).
\end{itemize}
Moreover, we construct explicit examples demonstrating that these bounds are attained.

We first address the case of quadratic first integrals on \(\mathbb{T}^2\). The following theorem characterizes completely the existence and uniqueness of simple crossing limit cycles of types \(\emph{\textbf{aa}}\) and \(\emph{\textbf{bb}}\).
\begin{theorem}\label{teorema1}
Let $X$ be a nonzero piecewise‑smooth vector field on the torus $\mathbb{T}^2$, given in $[0,1]^2$ by
\[
X(x,y)=
\begin{cases}
X_H(x,y), & (x,y)\in(0,1)^2,\\[4pt]
X_H\bigl(\pi(x,y)\bigr), & (x,y)\in\partial[0,1]^2,
\end{cases}
\]
and let $H(x,y) = ax^2 + bxy + cy^2$
be a quadratic first integral,
with discontinuity region exactly the boundary of the square  $\Sigma = \partial[0,1]^2$, and \(\pi\colon\partial[0,1]^2\to[0,1]^2\) is the usual projection that glues the edges. Then:

\begin{enumerate}
    \item $X$ can have at most one crossing limit cycles of type $\emph{\textbf{aa}}$ or $\emph{\textbf{bb}}$.
    \item A limit cycle of type $\emph{\textbf{bb}}$ exists if and only if the parameters $a$, $b$, and $c$ satisfy the following conditions:
    \begin{enumerate}
        \item $b \neq 0$, $c\neq 0$, $bc < 0$, and $|b|>|c|$.
        \item $\dfrac{(a c - b^2) (a + b - c)}{b} \leq 0$
        \item $ac<0$ or $ac>b^2$.
        \item If $\Delta \ge 0$, then $a(b - c)(b^2 + ab - ac) \ge 0$, where $\Delta = {b}^{4}-4\,a{b}^{2}c+4\,a{c}^{3}.$

    \end{enumerate}
\end{enumerate}
\end{theorem}

The previous result raises the question about first integrals of arbitrary degree. Let
\[
H(x,y)=\sum_{k=1}^n\sum_{j=0}^k a_{k,j}\,x^{k-j}y^j,
\]
and consider the same vector field defined on \((0,1)^2\). The next theorem gives an upper bound for the number of \(\emph{\textbf{bb}}\)-type cycles in terms of the degree \(n\).
\begin{theorem} \label{teo1}
Let $X$ be a nonzero piecewise-smooth vector field on the torus $\mathbb{T}^2$, defined on the fundamental domain $[0,1]^2$ by
\[
X(x,y)=
\begin{cases}
X_H(x,y), & (x,y)\in(0,1)^2,\\[4pt]
X_H\bigl(\pi(x,y)\bigr), & (x,y)\in\partial[0,1]^2,
\end{cases}
\]
with discontinuity region exactly the boundary of the square $\Sigma = \partial[0,1]^2$, and \(\pi\colon\partial[0,1]^2\to[0,1]^2\) is the projection that glues the edges . Let $H(x,y) = \sum_{k=1}^{n}\sum_{j=0}^{k}a_{k,j}\,x^{k-j}\,y^{j}$ be a polynomial first integral of degree $n \geq 1$. Assume that:
\begin{enumerate}
    \item $H(0,0) = 0$ and $H(x,0) - H(x,1) \not\equiv 0$
    \item For each $\emph{\textbf{bb}}$-type cycle $\gamma_{x_0} \subset \{H=H(x_0,0)\}$, $\min_{(x,y)\in\gamma_{x_0}} \|\nabla H(x,y)\| > 0$.

    \item For each such cycle, we have the conditions
          \[
          X(x_0,0) \cdot (0,1) \neq 0 \quad \text{and} \quad X(x_0,1) \cdot (0,-1) \neq 0.
          \]
\end{enumerate}
Then $X$ can have at most $(n-1)$ limit cycles of type $\emph{\textbf{bb}}$ (and analogously for type $\emph{\textbf{aa}}$).
\end{theorem}

While Theorem~\ref{teo1} covers cycles of type $\emph{\textbf{bb}}$ for general degree~$n$, quadratic integrals ($n=2$) exhibit richer looping patterns. The next result details bounds and existence conditions for $\emph{\textbf{aba}}$ and $\emph{\textbf{bab}}$ cycles when $H$ is quadratic.

\begin{theorem}\label{teo3}
Let $X$ be a nonzero piecewise‑smooth vector field on the torus $\mathbb{T}^2$, given in $[0,1]^2$ by
\[
X(x,y)=
\begin{cases}
X_H(x,y), & (x,y)\in(0,1)^2,\\[4pt]
X_H\bigl(\pi(x,y)\bigr), & (x,y)\in\partial[0,1]^2,
\end{cases}
\]
and let $H(x,y) = ax^2 + bxy + cy^2$
be a quadratic first integral,
with discontinuity region exactly the boundary of the square  $\Sigma = \partial[0,1]^2$, and \(\pi\colon\partial[0,1]^2\to[0,1]^2\) is the projection that glues the edges. Then:
\begin{enumerate}
    \item The vector field \( X \) can have at most two crossing limit cycles of types \( \emph{\textbf{aba}} \) or \( \emph{\textbf{bab}} \).
    
    \item A limit cycle of type \( \emph{\textbf{aba}} \) exists if and only if the parameters \( a \), \( b \), and \( c \) satisfy the following conditions:

\subsection*{Case 1: \protect\( a < 0 \protect\)}

    In this case, the conditions on  b and c are as follows:

\begin{enumerate}
\item If \( 2a < b < -\frac{2}{3} \sqrt{\frac{2}{3}} |a| \), then
        \[
        -\frac{1}{2} (a + b) + \frac{1}{2} \sqrt{a^2 + 2ab + 5b^2} < c \leq \rho_1. 
        \]
        \item If \( -\frac{2}{3} \sqrt{\frac{2}{3}} |a| \leq b < 0 \), then
        \[
        -\frac{1}{2} (a + b) + \frac{1}{2} \sqrt{a^2 + 2ab + 5b^2} < c \leq \rho_3.
        \]

        \item If \( 0 < b < -\frac{a}{2} \), then
        \[
        \rho_2 \leq c < \frac{-a^2 - ab + b^2}{a}.
        \]
    \end{enumerate}

    \subsection*{Case 2: \( a > 0 \)}

   In this case, the conditions on \( b \) and \( c \) are as follows:
    \begin{enumerate}
        \item If \( -\frac{a}{2} < b < 0 \), then
        \[
        \frac{-a^2 - ab + b^2}{a} < c \leq \rho_2.
        \]

        \item If \( 0 < b < 2a \), then
        \[
        \rho_1 \leq c < \frac{1}{2} (-a - b) - \frac{1}{2} \sqrt{a^2 + 2ab + 5b^2}.
        \]
    \end{enumerate}
\end{enumerate}
\noindent where \( P(c) \) is the polynomial defined by
\[
P(c) =  4a c^3 + 8a^2 c^2 + (4a^3 - 8ab^2) c + b^4,
\]
and \(\rho_k\) denotes the \( k \)-th root of the polynomial \( P(c) \), i.e., the value of \( c \) such that \( P(c) = 0 \), corresponding to the \( k \)-th solution of the equation \( P(c) = 0 \).
\end{theorem}

Finally, for arbitrary degree~$n$, we state an upper bound on the number of $\emph{\textbf{aba}}$-type cycles without requiring explicit existence criteria:

\begin{theorem} \label{teo2}
Let $X$ be a nonzero piecewise‑smooth vector field on the torus $\mathbb{T}^2$, given in $[0,1]^2$ by
\[
X(x,y)=
\begin{cases}
X_H(x,y), & (x,y)\in(0,1)^2,\\[4pt]
X_H\bigl(\pi(x,y)\bigr), & (x,y)\in\partial[0,1]^2,
\end{cases}
\]
with discontinuity region $\Sigma = \partial[0,1]^2$, and \(\pi\colon\partial[0,1]^2\to[0,1]^2\) is the projection that glues the edges. Let $H(x,y) = \sum_{k=1}^{n} \sum_{j=0}^{k} a_{k,j}  x^{k-j}  y^{j}$ be a polynomial first integral of degree $n$. Assume:
\begin{enumerate}
    \item $H(0,0) = 0$.
    \item For each $\emph{\textbf{aba}}$-type limit cycle $\gamma$ with transition points $p_1 = (0,y_1)$, $p_2 = (x_1,1)$, $p_3 = (x_1,0)$, $p_4 = (1,y_1)$:
    \begin{itemize}
        \item $\min_{(x,y)\in\gamma} \|\nabla H(x,y)\| > 0$,
         \item Transversality at sewing points:
        \[
        \langle X(p_i), n_{p_i} \rangle \neq 0 \quad \text{for } i=1,2,3,4,
        \]
        where $n_{p_i}$ is the inward-pointing unit normal to $\Sigma = \partial[0,1]^2$ at $p_i$.
    \end{itemize}
\end{enumerate}
Then $X$ admits at most $n(n-1)$ limit cycles of type $\emph{\textbf{aba}}$.
\end{theorem}

In the previous theorem, we established an upper bound for the number of limit cycles of type \( \emph{\textbf{aba}} \) in vector fields derived from first integrals of degree \( n \). A richer variety of cycle configurations can emerge when trajectories complete multiple vertical loops before returning to the initial position. We will not discuss these cases in the present paper.

The remainder of this paper is organized as follows. In Section 2 we enunciate the necessary preliminar results, including the basics about the theory of planar piecewise–smooth vector fields, first integrals, closing equations, and the homogeneization in \(\C P^2\) together with Bézout’s theorem. Section 3 is devoted to the proofs of our main results (Theorems \ref{teorema1}–\ref{teo2}), establishing upper bounds and existence criteria for various types of limit cycles on the torus. Finally, Section 4 presents explicit examples that illustrate and achieve the bounds provided by Theorems \ref{teo1} and \ref{teo2}.

\section{Preliminaries}
In this section we gather fundamental concepts and results that will be used throughout the paper.

\subsection{Planar Piecewise–Smooth Vector Fields}

Let \(\Sigma\subset\mathbb{R}^2\) be a smooth curve, defined as the zero set of a regular function \(h\colon\mathbb{R}^2\to\mathbb{R}\) with \(\nabla h(p)\neq0\) for all \(p\in\Sigma\). A \emph{planar piecewise–smooth vector field} is a mapping  
\[
X\colon\mathbb{R}^2\longrightarrow \mathbb{R}^2
\]
specified by two smooth vector fields \(X^+\) and \(X^-\) as follows:
\[
X(p)=
\begin{cases}
X^+(p), & h(p)>0,\\
X^-(p), & h(p)<0.
\end{cases}
\]
Solutions follow \(\dot z=X^\pm(z)\) in the regions \(\{h>0\}\) and \(\{h<0\}\). On the discontinuity curve \(\Sigma\), Filippov’s convention is applied to determine trajectory behavior.

Under this convention, \(\Sigma\) is partitioned into three disjoint subsets based on the signs of \(X^+h(p)\) and \(X^-h(p)\) relative to the normal gradient \(n=\nabla h\):

\begin{itemize}
  \item The \textbf{sewing region} (\emph{sewing region}) \(\Sigma^c\) is $\Sigma^c = \{p\in\Sigma : X^+h(p)\cdot X^-h(p)>0\}.$

  \item The \textbf{sliding region} (\emph{sliding region}) \(\Sigma^s\) is $\Sigma^s = \{p\in\Sigma : X^+h(p)<0 \text{ and } X^-h(p)>0\}.$

  \item The \textbf{escape region} (\emph{escape region}) \(\Sigma^e\) is $\Sigma^e = \{p\in\Sigma : X^+h(p)>0 \text{ and } X^-h(p)<0\}.$
\end{itemize}

This decomposition is essential for the qualitative analysis of the system, as it determines the possible orbit behaviors near the discontinuity.

In this work, we consider piecewise–smooth vector fields on the two-dimensional torus \(\mathbb{T}^2\), obtained by identifying opposite edges of the unit square \([0,1]^2\) via a projection \(\pi\). The vector field is given by
\[
X(x,y)=
\begin{cases}
X_H(x,y), & (x,y)\in(0,1)^2,\\[4pt]
X_H\bigl(\pi(x,y)\bigr), & (x,y)\in\partial[0,1]^2,
\end{cases}
\]
where \(H\colon[0,1]^2\to\mathbb{R}\) is a smooth first integral and \(\pi\colon\partial[0,1]^2\to[0,1]^2\) is the projection that glues the edges. Thus, the notions of sewing, sliding, and escape regions extend naturally to the discontinuities along the switching manifold \(\Sigma=\partial[0,1]^2\), analyzed locally as in the planar case and then interpreted globally on \(\mathbb{T}^2\).

\subsection{First integral}
Consider a differential system 
\begin{align}\label{int_prim_def}
\dot{x} = f(x), \hspace{1cm}x \in D \subset \mathbb{R}^n
\end{align}
where $D$ is an open subset of $\mathbb{R}^n$, and $f(x)=\left(f_1(x), f_2(x), \cdots , f_n(x)\right)$ is a $C^k$ function defined in $D$ with $k\geq 1$. A first integral of the differential system (\ref{int_prim_def}) is a continuous function $H(x)$ defined in the domain of definition $D$ of the differential system, which is not constant in any neighborhood, but such that $H(x)$ is constant on each orbit of the differential system (\ref{int_prim_def}). If $H(x)$ is $C^1$, then it is a first integral if and only if it satisfies:
\begin{align}
\dfrac{\partial H}{\partial x_1}f_1(x)+\dots+\dfrac{\partial H}{\partial x_n}f_n(x)=0,
\end{align}
for all the points $x \in D$.

\subsection{Closing equations}

\emph{Closing equations} are supplementary conditions that ensure continuity, uniqueness, or global consistency in piecewise-defined differential or algebraic systems, especially when discontinuities or regions with distinct evolution laws are present. In many applied contexts, these equations arise from physical principles, geometric criteria, or regularity requirements and are essential to coherently connect local regimes.

Consider a piecewise-smooth vector field \( Z = (X, Y) \) on \(\mathbb{R}^2\), with a discontinuity manifold \(\Sigma\) separating regions endowed with first integrals \(H_1(x,y)\) and \(H_2(x,y)\), associated with the vector fields \(X\) and \(Y\), respectively. The closing equation related to a periodic orbit that crosses \(\Sigma\) through a sequence of points
\[
\{(x_1,y_1),\,(x_2,y_2)\} \subset \Sigma
\]
is defined using the following conditions:
\begin{equation}\label{eq:fechamento}
\begin{aligned}
H_1(x_1,y_1) &= H_1(x_2,y_2),\\
H_2(x_1,y_1) &= H_2(x_2,y_2).\\
\end{aligned}
\end{equation}
This system expresses the continuity of the first integral along the cycle and is a necessary condition for the existence of periodic solutions that repeatedly cross the discontinuity. Together with the solution of these equations, we can check the direction of the trajectories, to assure that the intersection of the periodic orbit and $\Sigma$ occurs in a crossing region of $\Sigma$.

In the context of the torus, closing equations take an analogous form but are now related to the identifications of the fundamental square \( [0,1]^2 \). For instance, a closing condition arises when two points on opposite edges of the square, such as \( (x_1, 0) \) and \( (x_1, 1) \), are connected by a flow trajectory that preserves the level sets of a first integral \( H \). Since \( H \) remains constant along trajectories, a necessary condition for the existence of a limit cycle of type \( \emph{\textbf{bb}} \) is:
\[
H(x_1, 0) = H(x_1, 1).
\]
Similarly, for cycles of type \( \emph{\textbf{aba}} \), involving both horizontal and vertical crossings, the closing may involve four points, such as \( (0, y_1) \), \( (x_2, 0) \), \( (x_2, 1) \), and \( (1, y_1) \), leading to the system:
\[
\begin{aligned}
H(0, y_1) &= H(x_2, 1), \\
H(x_2, 0) &= H(1, y_1).
\end{aligned}
\]
The existence of real solutions to such equations is thus a sufficient condition for the presence of a limit cycle on the torus, with entry and exit points located on the discontinuity set \( \partial[0,1]^2 \).

\subsection{Homogenization in \(\C P^2\) and Bézout’s Theorem}
To count intersections of algebraic curves arising from the \emph{closing equations}, we employ techniques from the complex projective plane. Given an affine polynomial \(P(x,y)\) of degree \(d\), its \emph{homogenization} is defined by
\[
P^h(X,Y,Z) = Z^d \, P\!\Bigl(\tfrac X Z, \tfrac Y Z\Bigr),
\]
which is homogeneous of degree \(d\) in homogeneous coordinates \((X:Y:Z)\in\C P^2\). Thus each affine curve \(P=0\) corresponds to the projective curve \(P^h=0\).

We now state Bézout’s theorem in this projective setting:

\begin{bezout}[Bézout, {\cite[Thm.~I.3.1]{griffiths1994principles}}]
Let $C_1$ and $C_2$ be projective curves in $\C P^2$, defined by homogeneous polynomials of degrees $d_1$ and $d_2$, respectively, with no common component. Then $C_1$ and $C_2$ intersect in exactly $d_1\cdot d_2$ points in $\C P^2$, counted with multiplicity.
\end{bezout}

By ensuring there are no intersections at infinity and no shared components, we can apply Bézout directly to obtain sharp upper bounds on the total number of complex intersections. Subsequent hypotheses of smoothness of \(H\) and transversality along the discontinuity guarantee that all relevant intersections are real, simple, and lie strictly within the interior of the torus under consideration.

\section{Proof of the main results}
\subsection{Proof of Theorem \ref{teorema1}}

\begin{proof}[Proof of item (i)]
We show that for each crossing type there is at most one orbit from one side of the square to the opposite side staying in \((0,1)\times(0,1)\).
 
A cycle of type \(\emph{\textbf{bb}}\) must run from \((x_0,0)\) to \((x_0,1)\) with \(x_0\in(0,1)\). By invariance of \(H\):
\[
H(x_0,0) = H(x_0,1)
\;\Longrightarrow\;
a x_0^2 = a x_0^2 + b x_0 + c
\;\Longrightarrow\;
b x_0 + c = 0
\;\Longrightarrow\;
x_0 = -\frac{c}{b}.
\]
Thus at most one \(x_0\in(0,1)\) satisfies this (when \(b\neq0\) and \(-c/b\in(0,1)\)); otherwise no type-\(\emph{\textbf{bb}}\) cycle exists. Hence at most one.

The argument is analogous: a cycle must run from \((0,y_0)\) to \((1,y_0)\), giving
\[
H(0,y_0) = H(1,y_0)
\;\Longrightarrow\;
c y_0^2 = a + b y_0 + c y_0^2
\;\Longrightarrow\;
a + b y_0 = 0
\;\Longrightarrow\;
y_0 = -\frac{a}{b},
\]
so at most one \(y_0\in(0,1)\) satisfies this. 
 
Each seam‐closing equation is linear, yielding at most one solution in \((0,1)\). Therefore \(X_H\) has at most one crossing limit cycle of type \(\emph{\textbf{aa}}\) and at most one of type \(\emph{\textbf{bb}}\).
\end{proof}

\begin{proof}[Proof of item (ii)]
We already know from item (i) that, for limit cycles of type \emph{\textbf{bb}}, there is at most one candidate obtained by the “sewing closing” at
\[
x_0 = -\frac{c}{b}.
\]
Let us define the associated level constant:
\[
k = H(x_0,0) = a\,\frac{c^2}{b^2}.
\]
We will show that this level produces a limit cycle which crosses only at \((x_0,0)\) and \((x_0,1)\), without touching any other part of the boundary, precisely when conditions (a)–(d) are satisfied.

For the crossing at \(x_0\) to lie inside the interval \((0,1)\), it is necessary and sufficient that
\[
b \neq 0,
\quad
-\frac{c}{b}\in(0,1).
\]
This is equivalent to:
\[
b\neq0,
\quad
c\neq0,
\quad
bc<0,
\quad
|b|>|c|,
\]
which coincides with (a).

On the edge \(y=1\), the level curve satisfies:
\[
H(x,1) = ax^2 + bx + c = k
\quad\Longleftrightarrow\quad
P(x) = ax^2 + bx + c - \frac{ac^2}{b^2} = 0.
\]
We want the polynomial \(P(x)\) to have exactly one simple real root in the interval \((0,1)\), namely \(x_0 = -\tfrac{c}{b}\), and for no other root (if any) to lie in this interval. For this to occur, it is necessary and sufficient that
\[
\frac{(ac - b^2)\,(a + b - c)}{b} \le 0,
\]
which corresponds to condition (b).

At \(x=0\), we have:
\[
H(0,y) = c\,y^2 = k
\quad\Longleftrightarrow\quad
y^2 = \frac{ac}{b^2}.
\]
In order for there to be no solution \(y\in(0,1)\), we need:
\[
\frac{ac}{b^2} < 0
\quad\text{or}\quad
\frac{ac}{b^2} > 1
\quad\Longleftrightarrow\quad ac<0 \text{ or } ac>b^2,
\]
which is exactly (c).

On the edge \(x=1\), we consider:
\[
H(1,y) = c\,y^2 + b\,y + a = k
\quad\Longleftrightarrow\quad
Q(y) = c\,y^2 + b\,y + a - \frac{ac^2}{b^2} = 0.
\]
If the discriminant
\(\Delta = b^4 - 4a b^2 c + 4a c^3 < 0\), then there are no real roots and no intersection occurs.\\
If \(\Delta \ge 0\), the roots exist, and for both to lie outside \((0,1)\) it is necessary and sufficient that:
\[
Q(0)\,Q(1) > 0,
\]
which translates into the algebraic condition:
\[
a\,(b - c)\,(b^2 + ab - ac) > 0,
\]
that is, condition (d).

Since (a)–(d) guarantee that the level \(k\) crosses only at \((x_0,0)\) and \((x_0,1)\) and does not touch any other part of \(\partial[0,1]^2\), and there is a unique candidate (item (i)), we conclude that there exists a unique limit cycle of type \emph{\textbf{bb}} if and only if \((a,b,c)\) satisfy the stated conditions.
\end{proof}

\subsection{Proof of Theorem \ref{teo1}} \label{provateo1}
\begin{proof}
Consider a $\emph{\textbf{bb}}$-type limit cycle $\gamma_{x_0}$ crossing $y=0$ at $(x_0,0)$ and $y=1$ at $(x_0,1)$. Since $H$ is a first integral, it must be constant along trajectories, implying the necessary periodicity condition $H(x_0,0) = H(x_0,1)$. This equality defines the discriminant polynomial:
\[
P(x) = H(x,0) - H(x,1) = \sum_{k=1}^{n} \sum_{j=1}^{k} a_{k,j} x^{k-j}
\]
which has degree at most $n-1$ because the leading $x^n$ terms cancel in the difference. 

For each root $x_0 \in (0,1)$ of $P(x) = 0$, the non-degeneracy condition $\min_{(x,y)\in\gamma_{x_0}} \|\nabla H(x,y)\| > 0$ ensures $\gamma_{x_0}$ lies in a regular level set without critical points. Crucially, the transversality conditions
\[
X(x_0,0) \cdot (0,1) \neq 0 \quad \text{and} \quad X(x_0,1) \cdot (0,-1) \neq 0
\]
guarantee Filippov's crossing mode at both boundaries $y=0$ and $y=1$, ensuring trajectories properly cross the discontinuity rather than exhibiting sliding behavior or tangency.

By the Fundamental Theorem of Algebra, $P(x)$ admits at most $n-1$ real roots. Roots at $x=0$ or $x=1$ are excluded due to both the vertex identifications in the torus topology and potential singularities at these points. Thus, only roots in $(0,1)$ satisfying the non-degeneracy and transversality conditions yield valid $\emph{\textbf{bb}}$-cycles, bounding their number by $n-1$.The case $P(x) \equiv 0$ is handled by noting:
\begin{itemize}
    \item if $\nabla H \equiv 0$ on some $\gamma_{x_0}$, the non-degeneracy condition fails;
    \item if $\nabla H \not\equiv 0$, the system admits a continuum of periodic orbits, which are necessarily non-isolated.
\end{itemize}
The analogous result for $\emph{\textbf{aa}}$-cycles follows by considering $Q(y) = H(0,y) - H(1,y)$.
\end{proof}

\subsection{Proof of Theorem \ref{teo3}}
To bound the number of limit cycles, we impose periodicity conditions on the torus \(\mathds{T}^2\).  Consider the closing equations
\begin{align}
\begin{cases}\label{equafechamento}
H(0, y) = H(x, 1),\\
H(x, 0) = H(1, y).
\end{cases}
\end{align}

Algebraic resolution of \eqref{equafechamento} yields two nontrivial solutions:

\begin{align}\label{solu1}
\begin{cases} 
x_1 = \dfrac{-b^2 + 2c(a + c) + Q}{2b(a - c)}, \\[8pt]
y_1 = \dfrac{b^2 - 2a(a + c) - Q}{2b(a - c)},
\end{cases}
\qquad
\begin{cases}
x_2 = \dfrac{-b^2 + 2c(a + c) - Q}{2b(a - c)}, \\[8pt]
y_2 = \dfrac{b^2 - 2a(a + c) + Q}{2b(a - c)},
\end{cases}
\end{align}
where 
\[
Q \;=\;\sqrt{\,b^4 \;-\;8ab^2c\;+\;4ac(a+c)^2\,},
\]
governs the existence of limit cycles: if \(Q>0\), there are two distinct real cycles (of types \(\emph{\textbf{aba}}\) and \(\emph{\textbf{bab}}\)); if \(Q=0\), a degenerate cycle occurs; and if \(Q<0\), no limit cycle is present.
\begin{remark} \label{obs1}
Since the torus is obtained by identifying the opposite edges of the unit square \([0,1]^2\), the intersection points \((x_i,y_i)\) given in \eqref{solu1} must lie strictly in the interior of that square in order to represent well-defined points on the torus:
\[
(x_i,y_i)\in(0,1)^2,\quad i=1,2.
\]
\end{remark}

Note also that
\[
Q \;=\;\sqrt{b^4 - 8ab^2c + 4ac(a+c)^2}
\]
is the discriminant governing the existence of limit cycles: when \(Q>0\) there are two distinct real cycles (of types \emph{\emph{\textbf{aba}}} and \emph{\textbf{bab}}); if \(Q=0\), the system exhibits a degenerate cycle; and if \(Q<0\), no limit cycles occur.

From this expression we obtain the solutions in \eqref{solu1} and the cubic polynomial
\[
  P(c)\;=\;Q^2 \;=\;4a\,c^3 + 8a^2\,c^2 + (4a^3-8ab^2)\,c + b^4,
\]
whose ordered real roots \(\rho_1<\rho_2<\rho_3\) determine the sign of \(P(c)\).  To analyze \(P(c)\), we study its discriminant and the nature of its roots.  The discriminant is given by  
\[
\Delta(P) \;=\; 16\,a^2\,b^2\,(2a - b)^2\,(2a + b)^2\,(8a^2 - 27b^2)\,.
\]

If \(\Delta(P)>0\), the polynomial has three distinct real roots; if \(\Delta(P)=0\), it has multiple roots (either a triple root or one double root and one simple root); and if \(\Delta(P)<0\), it has one real root and two complex conjugate roots.

Before treating each subcase according to the sign of \(\Delta(P)\), it is useful to recall how the leading term of \(P(c)\) dominates the global behavior of the polynomial.  The following lemma summarizes this property, allowing us to determine immediately the sign of \(P(c)\) as \(c\to\pm\infty\) based solely on the sign of \(a\).

\begin{lema}\label{lem:asymp}
Let
\[
P(c) \;=\; 4a\,c^3 \;+\; 8a^2c^2 \;+\;(4a^3 - 8ab^2)\,c \;+\; b^4.
\]
Then
\[
\lim_{c\to +\infty} P(c)
=
\begin{cases}
-\infty, & a<0,\\[1mm]
+\infty, & a>0,
\end{cases}
\qquad
\lim_{c\to -\infty} P(c)
=
\begin{cases}
+\infty, & a<0,\\[1mm]
-\infty, & a>0.
\end{cases}
\]
In particular, the sign of \(P(c)\) for large \(\lvert c\rvert\) is determined solely by the sign of \(a\).
\end{lema}
In all of the following cases, we invoke Lemma~\ref{lem:asymp} together with the continuity of \(P(c)\) to identify exactly those intervals of \(c\) for which \(P(c)\ge0\).  We then intersect these intervals with the geometric requirement \((x_i,y_i)\in(0,1)^2\) to single out the parameter values that yield admissible limit‐cycle solutions.

\subsection{Case \(\Delta(P)<0\)}
\label{sec:DeltaNeg}\leavevmode  

Assume that \(\Delta(P)<0\), so that the cubic polynomial \(P(c)\) has exactly one real root.  In this regime there are two possibilities for the sign of \(a\), each giving rise to four disjoint intervals for \(b\):

\begin{align*}
&\text{Subcase 1.1: } a<0,\quad b\in(-\infty,2a)\,\cup\,(2a,-\tfrac{2}{3}\sqrt{\tfrac{2}{3}}|a|)\,\cup\,(\tfrac{2}{3}\sqrt{\tfrac{2}{3}}|a|,-2a)\,\cup\,(-2a,\infty),\\
&\text{Subcase 1.2: } a>0,\quad b\in(-\infty,-2a)\,\cup\,(-2a,-\tfrac{2}{3}\sqrt{\tfrac{2}{3}}a)\,\cup\,(\tfrac{2}{3}\sqrt{\tfrac{2}{3}}a,2a)\,\cup\,(2a,\infty).
\end{align*}

\subsubsection{Subcase \(a<0\)}

In Subcase 1.1 (\(a<0\)), Lemma~\ref{lem:asymp} together with continuity of \(P(c)\) imply that its unique real root $\rho_1 $ satisfies
\[
P(c)\ge0
\qquad\text{if and only if}\qquad
c\in(-\infty,\rho_1].
\]

Moreover, since the torus \(\mathds T^2\) is obtained by identifying opposite edges of the unit square \((0,1)^2\), any solution \((x_i,y_i)\) of \eqref{solu1} must lie strictly in the interior of that square to define a genuine point on the torus (see Remark ~\ref{obs1}).

We therefore must impose further restrictions on the parameters so that both \((x_1,y_1)\) and \((x_2,y_2)\) remain in \((0,1)^2\).  We now examine each of the four intervals for \(b\).

\subsubsection{Interval \(b\in(-\infty,2a)\)}

For \(x_1\in(0,1)\) one needs
\[
-a+b < c < \frac{-a-b-\sqrt{a^2+2ab+5b^2}}{2}
\quad\text{or}\quad
-a-b < c < \frac{-a-b+\sqrt{a^2+2ab+5b^2}}{2},
\]
while \(x_2\in(0,1)\) requires
\[
0 < c < -a-b.
\]
Since these two sets of \(c\) do not overlap, no choice of \(c\) in this range can make both \(x_1\) and \(x_2\) lie in \((0,1)\).  Thus \(b\in(-\infty,2a)\) yields no admissible torus solutions.

\subsubsection{Subcase \(b \in \bigl(2a,\,-\tfrac{2}{3}\sqrt{\tfrac{2}{3}}\,|a|\bigr)\)}\label{subcaso(a1)}

For \(x_1\in(0,1)\), one must have
\[
-a - b < c \le \rho_1.
\]
For \(x_2\in(0,1)\), one of the following must hold:
\[
\begin{aligned}
\frac{-a - b - \sqrt{a^2 + 2ab + 5b^2}}{2} &< c < 0,\\
-a + b &< c < -a - b,\\
\frac{-a - b + \sqrt{a^2 + 2ab + 5b^2}}{2} &< c \le \rho_1.
\end{aligned}
\]
Analogously, \(y_1\in(0,1)\) requires
\[
-a - b < c \le \rho_1,
\]
while \(y_2\in(0,1)\) holds if and only if
\[
-a + b < c < \frac{-a^2 - ab + b^2}{a}
\quad\text{or}\quad
-a - b < c \le \rho_1.
\]
Intersecting these four families of inequalities shows that, for
\[
b\in\bigl(2a,\,-\tfrac{2}{3}\sqrt{\tfrac{2}{3}}\,|a|\bigr),
\]
the only admissible \(c\) satisfy
\[
\frac{-a - b + \sqrt{a^2 + 2ab + 5b^2}}{2}
\;<\; c \;\le\; \rho_1.
\]
For the complementary ranges \(b\in\bigl(\tfrac{2}{3}\sqrt{\tfrac{2}{3}}\,|a|,\,-2a\bigr)\) or \(b\in(-2a,\infty)\), no single choice of \(c\) meets all four conditions \(x_i,y_i\in(0,1)\), so there are no valid solutions on the torus in those regions.

In summary, a valid torus solution exists precisely when
\[
a<0,\quad 2a < b < -\tfrac{2}{3}\sqrt{\tfrac{2}{3}}\,|a|,\quad
\tfrac{-a - b + \sqrt{a^2 + 2ab + 5b^2}}{2}
< c \le \rho_1,
\]
where \(\rho_1=\operatorname{Root}(P(c),1)\).  These conditions ensure (i) that \(P(c)\ge0\) and (ii) that both \(x_i\) and \(y_i\) lie strictly in \((0,1)\), completing the proof in this subcase.  

\subsubsection{Subcase \(a>0\)}

For \(a>0\), Lemma~\ref{lem:asymp} together with the continuity of \(P(c)\) implies that, if we denote by $\rho_2 \;=\;\operatorname{Root}\bigl(P(c),2\bigr)$
the smallest real number for which \(P(c)=0\) and \(P(c)\ge0\) thereafter, then
\[
P(c)\ge0
\quad\Longleftrightarrow\quad
c\in[\rho_2,\infty).
\]

\subsubsection{Case \(b \in \bigl(\tfrac{2}{3}\sqrt{\tfrac{2}{3}}\,a,\,2a\bigr)\)}\label{subcase(b3)}

In this subcase, we also require that the formulas in \eqref{solu1} yield points \((x_i,y_i)\in(0,1)^2\). Specifically, \(x_1\in(0,1)\) forces
\[
\rho_1\;\le\;c\;<\;-a-b,
\]
whereas \(x_2\in(0,1)\) holds under one of the following three alternatives:
\[
\begin{aligned}
\rho_1\;\le\;&c\;<\;\tfrac{-a-b}{2}-\tfrac12\sqrt{a^2+2ab+5b^2},\\
-a-b\;<\;&c\;<\;0,\\
-a+b\;<\;&c\;<\;\tfrac{-a-b}{2}+\tfrac12\sqrt{a^2+2ab+5b^2}.
\end{aligned}
\]
Similarly, \(y_1\in(0,1)\) imposes the same bound as \(x_1\), and \(y_2\in(0,1)\) requires
\[
\rho_1\;\le\;c\;<\;-a-b
\quad\text{or}\quad
\tfrac{-a^2-ab+b^2}{a}\;<\;c\;<\;-a+b.
\]
Intersecting all these conditions shows that, whenever
\[
\tfrac{2}{3}\sqrt{\tfrac{2}{3}}\,a < b < 2a,
\]
the unique interval of \(c\) satisfying both \(P(c)\ge0\) and \((x_i,y_i)\in(0,1)^2\) is
\[
\rho_1\;\le\;c\;<\;\tfrac{-a-b}{2}-\tfrac12\sqrt{a^2+2ab+5b^2}.
\]
For any other values of \(b\) in \(\bigl(-\infty,-2a\bigr)\), \(\bigl(-2a,\tfrac{2}{3}\sqrt{\tfrac{2}{3}}\,a\bigr)\), or \((2a,\infty)\), no \(c\) can satisfy all four inequalities simultaneously, ruling out valid torus solutions in those cases.

In summary, in the regime \(a>0\) and \(b\in\bigl(\tfrac{2}{3}\sqrt{\tfrac{2}{3}}\,a,2a\bigr)\), the necessary and sufficient conditions for solutions to lie in \((0,1)^2\) are
\[
a>0,\quad \tfrac{2}{3}\sqrt{\tfrac{2}{3}}\,a<b<2a,\quad
\rho_1\le c<\tfrac{-a-b}{2}-\tfrac12\sqrt{a^2+2ab+5b^2}.
\]

\subsection{Case \(\Delta(P)>0\)}
\label{sec:DeltaPos}\leavevmode  

Suppose now that \(\Delta(P)>0\), so that the cubic polynomial \(P(c)\) has three distinct real roots. In this regime, the admissible ranges of \(b\) are the complements of those in the case \(\Delta(P)<0\).  Thus we set
\[
\begin{aligned}
\text{Subcase 2.1: } a<0,\quad &b\in\bigl(-\tfrac{2}{3}\sqrt{\tfrac{2}{3}}|a|,0\bigr)\;\cup\;\bigl(0,\tfrac{2}{3}\sqrt{\tfrac{2}{3}}|a|\bigr),\\[1mm]
\text{Subcase 2.2: } a>0,\quad &b\in\bigl(-\tfrac{2}{3}\sqrt{\tfrac{2}{3}}\,a,0\bigr)\;\cup\;\bigl(0,\tfrac{2}{3}\sqrt{\tfrac{2}{3}}\,a\bigr).
\end{aligned}
\]
These intervals of \(b\) are exactly those for which \(\Delta(P)>0\) ensures three real roots of \(P(c)=0\) and hence the possibility of torus–valid solutions.

\subsection{Subcase \(a<0\)}

When \(a<0\) and \(\Delta(P)>0\), Lemma~\ref{lem:asymp} gives
\[
P(c)\to -\infty\quad(c\to+\infty),
\qquad
P(c)\to +\infty\quad(c\to-\infty).
\]
Let \(\rho_1<\rho_2<\rho_3\) be the three real zeros of \(P(c)\).  By continuity and sign–changes,
\[
P(c)\ge0
\quad\Longleftrightarrow\quad
c\in(-\infty,\rho_1]\;\cup\;[\rho_2,\rho_3].
\]
We then restrict to \(b\in\bigl(-\tfrac{2}{3}\sqrt{\tfrac{2}{3}}|a|,0\bigr)\cup\bigl(0,\tfrac{2}{3}\sqrt{\tfrac{2}{3}}|a|\bigr)\), choose \(c\) in these nonnegative regions, and simultaneously impose \((x_i,y_i)\in(0,1)^2\).  The intersection of these requirements yields the precise admissible bands for \(c\) in each subinterval of \(b\).

\subsubsection{Subcase \(b\in\bigl(-\tfrac{2}{3}\sqrt{\tfrac{2}{3}}|a|,0\bigr)\)}\label{subcaso(a2)}

One checks that
\[
  -a - b < c \le \rho_3
\]
is necessary for \(x_1\in(0,1)\), while \(x_2\in(0,1)\) forces one of
\[
  \begin{aligned}
    \frac{-a-b-\sqrt{a^2+2ab+5b^2}}{2} &< c < 0,\\
    -a+b                           &< c < -a-b,\\
    \frac{-a-b+\sqrt{a^2+2ab+5b^2}}{2} &< c \le \rho_3.
  \end{aligned}
\]
Analogous inequalities for \(y_1,y_2\) intersect to
\[
  \frac{-a-b+\sqrt{a^2+2ab+5b^2}}{2} < c \le \rho_3.
\]
Hence, for \(-\tfrac{2}{3}\sqrt{\tfrac{2}{3}}|a|<b<0\) one obtains
\[
  \frac{-a-b+\sqrt{a^2+2ab+5b^2}}{2} < c \le \rho_3.
\]

\subsubsection{Subcase \(b\in\bigl(0,\tfrac{2}{3}\sqrt{\tfrac{2}{3}}|a|\bigr)\)}\label{subcaso(a3)}
Similarly, enforcing \(x_1,x_2,y_1,y_2\in(0,1)\) yields several alternative bands whose intersection is exactly
\[
\rho_2\;\le\;c\;<\;\frac{-a^2-ab+b^2}{a},
\]
valid when \(0<b<-\tfrac{a}{2}\).  This completes the \(a<0\) analysis.

\subsection{Subcase \(a>0\)}

If \(a>0\) and \(\Delta(P)>0\), Lemma~\ref{lem:asymp} and continuity imply
\[
P(c)\ge0
\quad\Longleftrightarrow\quad
c\in[\rho_1,\rho_2]\;\cup\;[\rho_3,\infty).
\]
We then take \(b\in\bigl(-\tfrac{2}{3}\sqrt{\tfrac{2}{3}}\,a,0\bigr)\cup\bigl(0,\tfrac{2}{3}\sqrt{\tfrac{2}{3}}\,a\bigr)\), identify the regions of \(c\) with \(P(c)\ge0\), and intersect with the geometric constraints \(x_i,y_i\in(0,1)\). 

This combined algebraic–geometric argument yields the exact parameter bands supporting limit cycles of types \emph{\textbf{aba}} or \emph{\textbf{bab}} when \(\Delta(P)>0\).

\subsubsection{Subcase \(b \in \Bigl(-\frac{2}{3}\sqrt{\frac{2}{3}}\,a,0\Bigr)\)}\label{subcase(b1)}

To guarantee \(x_1 \in (0,1)\), it is essential to verify
\[
\begin{aligned}
-a+b &< c < \frac{-a-b}{2} - \frac{1}{2}\sqrt{a^2+2ab+5b^2}, \quad \text{or} \\
-a-b &< c \le \rho_2, \quad \text{or} \\
\rho_3 &\le c < \frac{-a-b}{2} + \frac{1}{2}\sqrt{a^2+2ab+5b^2}.
\end{aligned}
\]
For \(x_2 \in (0,1)\) to hold, one of the following conditions must occur:
\[
-a - b < c \le \rho_2 \quad \text{or} \quad \rho_3 \le c < 0.
\]
Similarly, to ensure \(y_1 \in (0,1)\), we require
\[
-a+b < c < -a-b \quad \text{or} \quad \frac{-a^2-ab+b^2}{a} < c \le \rho_3.
\]
whereas for \(y_2 \in (0,1)\) to be satisfied, it must hold that
\[
-a-b < c \le \rho_2.
\]
Additionally, for both \(y_1\) and \(y_2\) to lie in \((0,1)\), it is necessary that \(\frac{-a}{2}<b<0\). Taking the intersection of all these conditions, we conclude that for
\[
b \in \Bigl(-\frac{2}{3}\sqrt{\frac{2}{3}}\,a,0\Bigr)
\]
and for the solutions \((x_i,y_i)\) to remain within the torus, \(c\) must satisfy
\[
\frac{-a^2 - ab + b^2}{a} < c \le \rho_2.
\]
In summary, for valid intersection of solutions within the torus in this subcase, parameters must satisfy
\[
a < 0,\quad -\frac{a}{2} < b < 0,\quad \frac{-a^2 - ab + b^2}{a} < c \le \rho_2.
\]
These conditions ensure that variables \(x_i\) and \(y_i\) remain in \((0,1)\), preserving solution integrity within the torus.

\subsubsection{Subcase \(b \in \Bigl(0, \frac{2}{3}\sqrt{\frac{2}{3}}\,a\Bigr)\)}\label{subcase(b2)}

To satisfy \(x_1 \in (0,1)\), it is imperative that
\[
\rho_1 \le c < -a-b.
\]
For \(x_2 \in (0,1)\) to be guaranteed, one of the following must occur:
\[
\begin{aligned}
\rho_1 \le c &< \frac{-a-b}{2} - \frac{1}{2}\sqrt{a^2+2ab+5b^2}, \quad \text{or} \\
-a-b < c &< -a+b, \quad \text{or} \\
0 < c &< \frac{-a-b}{2} + \frac{1}{2}\sqrt{a^2+2ab+5b^2}.
\end{aligned}
\]
Similarly, to ensure \(y_1 \in (0,1)\), we require
\[
\rho_1 \le c < -a-b.
\]
whereas for \(y_2 \in (0,1)\) to hold, it must satisfy
\[
\rho_1 \le c < -a-b \quad \text{or} \quad \frac{-a^2-ab+b^2}{a} < c < -a+b.
\]
Taking the intersection of these conditions, we conclude that for \(0< b <\frac{2}{3}\sqrt{\frac{2}{3}}\,a\), and for solutions \((x_i,y_i)\) to remain within the torus, \(c\) must satisfy
\[
\rho_1 \le c < -\frac{a+b}{2} - \frac{1}{2}\sqrt{a^2+2ab+5b^2}.
\]
In summary, for valid solution intersections within the torus in this subcase, parameters must satisfy
\[
a > 0,\quad 0 < b < \frac{2}{3}\sqrt{\frac{2}{3}}\,a,\quad \rho_1\le c < -\frac{a+b}{2} - \frac{1}{2}\sqrt{a^2+2ab+5b^2}.
\]
These conditions ensure that solutions \((x_i,y_i)\) respect the \((0,1)\) interval constraint, maintaining result integrity within the torus.

\subsection{Case \(\Delta(P)=0\)}
\label{sec:DeltaZero} \leavevmode  

Let us now suppose that \(\Delta(P)=0\), i.e., the polynomial \(P(c)\) has at least one real root with multiplicity greater than one. By Lemma~\ref{lem:asymp}, its asymptotic behavior as \(c\to\pm\infty\) is entirely determined by the sign of \(a\). The presence of repeated roots strongly constrains the intervals where \(P(c)\ge0\), which in turn restricts the admissible values of \(c\) for solutions \((x_i,y_i)\) to remain in \((0,1)\). As demonstrated, critical cases for \(\Delta(P)=0\) occur at specific values of \(b\), distributed into two subcases according to the sign of \(a\):

\subsubsection{Subcase \(a<0\)}\label{d1}

\[
b \in \{\, 2a,\; \tfrac{2}{3}\sqrt{\tfrac{2}{3}}\,|a|,\; 0,\; -\tfrac{2}{3}\sqrt{\tfrac{2}{3}}\,|a|,\; -2a \,\}.
\]
It is observed that for \(b = 2a\), \(b = \tfrac{2}{3}\sqrt{\tfrac{2}{3}}\,|a|\), and \(b = -2a\), the intersection of solution sets is empty, and no solution exists for \(b=0\). For \(b = -\tfrac{2}{3}\sqrt{\tfrac{2}{3}}\,|a|\), the solution becomes
\[
\frac{-a-b+\sqrt{a^2+2ab+5b^2}}{2} < c \le \rho_3.
\]
This completes the analysis for \(a<0\).

\subsubsection{Subcase \(a>0\)}\label{d2}
\[
b \in \{\, -2a,\; -\tfrac{2}{3}\sqrt{\tfrac{2}{3}}\,a,\; 0,\; \tfrac{2}{3}\sqrt{\tfrac{2}{3}}\,a,\; 2a \,\}.
\]
It is verified that for \(b = -2a\), \(b = -\tfrac{2}{3}\sqrt{\tfrac{2}{3}}\,a\), and \(b = 2a\), the intersection of solution sets is empty, with no solution existing for \(b=0\). For \(b = \tfrac{2}{3}\sqrt{\tfrac{2}{3}}\,a\), the solution is given by
\[
\rho_1 \le c < \frac{-a-b}{2} - \frac{1}{2}\sqrt{a^2+2ab+5b^2}.
\]
This concludes the analysis for \(a>0\).

The overlap of solution sets, conditioned by these \(b\)-constraints and \(c\)-intervals determined from \(P(c)=0\) with multiple roots, ensures that only at these critical values can valid toroidal solutions exist. Thus, the \(\Delta(P)=0\) case analysis is complete.

The proof of Theorem~\ref{teo3} was structured into three distinct scenarios based on the discriminant \(\Delta(P)\). First, we showed that the torus-closing conditions (\ref{equafechamento}) permit at most two limit cycles (corresponding to \emph{\textbf{aba}} and \emph{\textbf{bab}} types), establishing item (1).  

Next, for the existence of an \emph{\textbf{aba}}-type cycle (item 2), we divided the analysis by the sign of \(a\). When \(a<0\), changes in the number and position of \(P(c)\)'s real roots led to three critical \(b\)-subcases: Subcase~\ref{subcaso(a1)} validates clause (a), Subcase~\ref{subcaso(a2)} leads to clause (b), and Subcase~\ref{subcaso(a3)} confirms clause (c). For \(a>0\), Subcase \ref{subcase(b1)} supports clause (a), while Subcases~\ref{subcase(b2)} and~\ref{subcase(b3)} jointly ensure clause (b). Throughout, Lemma~\ref{lem:asymp} and the continuity of \(P(c)\) enabled precise descriptions of \(c\)-intervals where \(P(c)\ge0\), with torus geometry requiring \((x_i,y_i) \in (0,1)^2\).  

Finally, the \(\Delta(P)=0\) study revealed critical \(b\)-values where the polynomial has multiple roots. At these boundaries, equalities in the \(c\)-constraints exactly match conditions examined in Subcases~\ref{d1} and~\ref{d2}, guaranteeing the degenerate cycle predicted in item 2.b.

Thus, we exhaustively covered all possible \((a,b,c)\) configurations, proving that the conditions stated in each theorem item are simultaneously necessary and sufficient for limit cycle existence in the \(X_H\) dynamics. The proof is complete.

\subsection{Proof of Theorem \ref{teo2}}

\begin{proof}
We denote by 
\[
P(x,y)=H(0,y)-H(x,1),
\qquad
Q(x,y)=H(x,0)-H(1,y)
\]
the two closing‐equation polynomials whose common real solutions in \((0,1)^2\) correspond bijectively to \(\emph{\textbf{aba}}\)–type limit cycles.  Since 
\[
H(x,y)=\sum_{k=1}^n\sum_{j=0}^k a_{k,j}\,x^{\,k-j}y^j
\]
has total degree \(n\), a inspection, $\deg_x P = \deg_x Q = n$, but with algebraic manipulations, the system $P = Q = 0$ is equivalent to $P = \widetilde{Q} = 0$, where $\deg_x \widetilde{Q} = n - 1$.  We now embed these affine curves into the complex projective plane \(\C P^2\) by homogenizing with respect to coordinates \((X:Y:Z)\):
\[
P^h(X,Y,Z)=Z^n\,P\!\Bigl(\frac XZ,\frac YZ\Bigr),
\qquad
Q^h(X,Y,Z)=Z^{\,n-1}\,\widetilde{Q}\!\Bigl(\frac XZ,\frac YZ\Bigr).
\]
By construction, \(P^h\) is homogeneous of degree \(n\) and \(Q^h\) of degree \(n-1\), and their common zeros in the affine chart \(Z\neq0\) coincide with the real solutions of \(P=\widetilde{Q}=0\).

To verify that no spurious intersections appear at infinity, we set \(Z=0\).  One finds
\[
P^h(X,Y,0)=a_{n,n}\,Y^n,
\qquad
Q^h(X,Y,0)=-\,a_{n,0}\,X^n.
\]
Since \(a_{n,0}\) and \(a_{n,n}\) cannot both vanish (else \(\deg H<n\)), these two homogeneous polynomials vanish on disjoint points of the line at infinity, so there is no common intersection there.

Now Bézout’s theorem in \(\C P^2\) ensures that two projective curves of degrees \(n\) and \(n-1\), having no common component and no intersection at infinity, meet in exactly \(n(n-1)\) points when counted with multiplicity.  Finally, the non‐degeneracy conditions 
\[
\min_{(x,y)\in\gamma}\|\nabla H(x,y)\|>0
\quad\text{and}\quad
\bigl\langle X(p_i),n_{p_i}\bigr\rangle\neq0
\quad(i=1,\dots,4)
\]
guarantee that all intersections in the affine chart are simple real points in \((0,1)^2\).  Hence there can be at most \(n(n-1)\) distinct real solutions of \(P=Q=0\), and therefore at most \(n(n-1)\) \(\emph{\textbf{aba}}\)–type limit cycles, as claimed.
\end{proof}

In the next section, we will present some numerical examples in which these upper bounds are reached. Theorem \ref{teo1} shows that this occurs for any value of $n \in \mathbb{N}$. However, for Theorem \ref{teo2}, we will see that the bounds are reached when $n$ equals 1 and 2.

\section{Examples}
In this section, we will present some examples demonstrating the validity of the bounds established in Theorems \ref{teo1} and \ref{teo2}.

Let's illustrate Theorem (\ref{teo2}) with the following example.

\begin{exemplo}\label{ex:degree3}
Consider the following nonlinear differential system \(X_H\in\mathfrak{X}^r(\mathbb{T}^2)\) on the torus, expressed in the fundamental square \([0,1]^2\):
{\footnotesize \begin{align}\label{sistema1}
 \left\{\begin{array}{l}
\dot{x}=0.002642824354\,{x}^{2}- 1.078842526\,xy+ 0.5370427213\,x+
0.2912880768\\
\\
\dot{y}=5.561516025\,{x}^{2}- 0.005285648708\,xy+ 0.5394212632\,{y}^{2}-
5.000078774\,x- 0.5370427213\,y+ 1.071920970.
\end{array}\right.
\end{align}}
One verifies that this system admits the cubic first integral
\begin{equation}\label{intsistema1}
\begin{split}
H(x,y) ={}& -1.853838675\,x^{3} + 0.002642824354\,x^{2}y - 0.5394212632\,x y^{2}
+ 2.500039387\,x^{2}\\
&\quad + 0.5370427213\,x y - 1.071920970\,x + 0.2912880768\,y.
\end{split}
\end{equation}
Since \(\deg H=3\), Theorem~\ref{teo2} predicts at most \(3\cdot2=6\) limit cycles of type \(\emph{\textbf{aba}}\).  We now show that this bound is attained.

Because trajectories lie on the level curves of \(H\), each \(\emph{\textbf{aba}}\)-type cycle corresponds to a solution \((x,y)\in(0,1)^2\) of the closing equations
\[
H(x,0)=H(1,y),
\qquad
H(0,y)=H(x,1).
\]
Numerical solution of these equations yields six distinct pairs, rounded to three decimal places:
\[
(x_1,y_1)\approx(0.25,\,0.516),\quad
(x_2,y_2)\approx(0.33,\,0.490),\quad
(x_3,y_3)\approx(0.41,\,0.494),
\]
\[
(x_4,y_4)\approx(0.49,\,0.507),\quad
(x_5,y_5)\approx(0.57,\,0.511),\quad
(x_6,y_6)\approx(0.65,\,0.485).\]
Each pair \((x_i,y_i)\) determines a unique \(\emph{\textbf{aba}}\)-cycle
\[
\gamma_i \;=\;\bigl\{H(x,y)=H(0,y_i)\bigr\}\,\cap\,\bigl\{H(x,y)=H(x_i,0)\bigr\},
\]
which is a smooth, simple closed curve on \(\mathbb{T}^2\).  Figure~\ref{todosjuntos} displays all six cycles simultaneously.

\begin{figure}[h]
  \centering
  \begin{overpic}[height=5cm,unit=1mm]{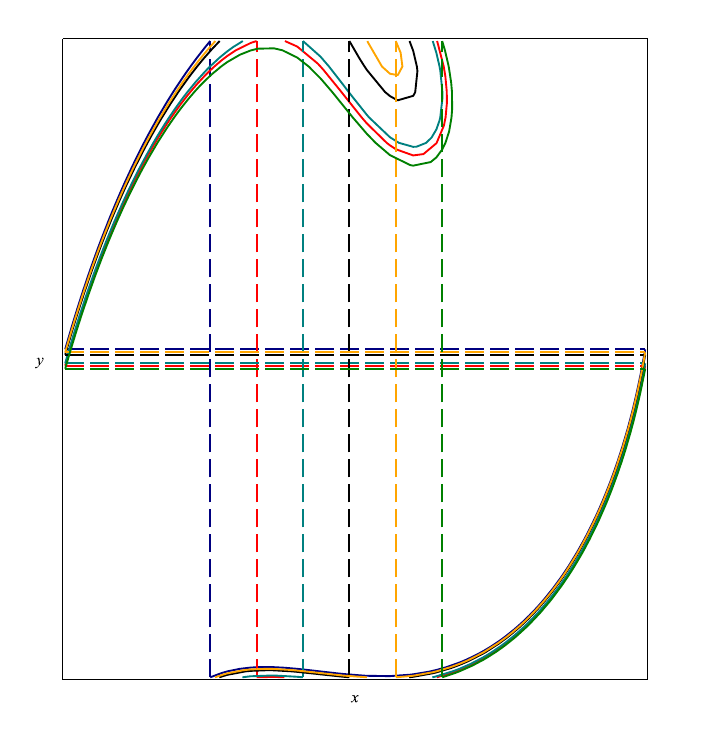}
    \put(43,20){\(\gamma_i\)}
  \end{overpic}
  \quad
  \begin{overpic}[height=5cm,unit=1mm]{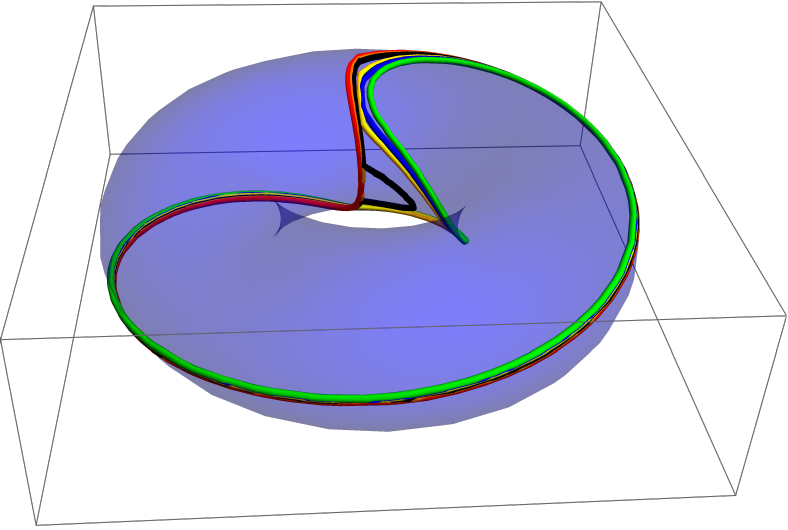}
    \put(25,28){\(\gamma_i\)}
  \end{overpic}
  \caption{All six \(\emph{\textbf{aba}}\)-type limit cycles of system \eqref{sistema1}.}
  \label{todosjuntos}
\end{figure}

Finally, one checks that each \(\gamma_i\) crosses the switching manifold transversely (i.e.\ it is a sewing cycle) by verifying
\[
\bigl\langle X_H(p),n(p)\bigr\rangle\neq0
\quad\text{at each intersection }p\in\partial[0,1]^2.
\]
Thus System \eqref{sistema1} realizes the maximum of six \(\emph{\textbf{aba}}\)-type limit cycles, confirming the sharpness of Theorem~\ref{teo2}.
\end{exemplo}

\begin{exemplo}

Consider the polynomial first integral of degree \(n\) given by
\[
H(x,y) = a_{0}\,x^{n} + y \prod_{i=1}^{n-1}\bigl(x - x_i\bigr),
\qquad
a_{0}=(-1)^{n+1},
\qquad
x_i=\frac{i}{n}.
\]
Let us take the associated Hamiltonian vector field,
\[
X_H = \bigl(-\partial_y H,\;\partial_x H\bigr),
\]
and observe that \(H\) is a first integral of \(X_H\). Without loss of generality, we may choose \(X_H\) to represent all vector fields sharing this same first integral, since any such vector field can be written as
\[
X = \mu(x,y)\,X_H,
\]
for some smooth scalar function \(\mu(x,y)\). In this way, all such dynamics share the same foliation given by the level curves of \(H\), differing only in the parametrization of time along the trajectories, which is governed by the function \(\mu\).

With the choice \(a_{0}=(-1)^{n+1}\), we note that \(H(0,0)=0\). For each \(k=1,\dots,n-1\), the vertical line \(x_k = k/n\) satisfies
\[
H(x_k,0)=H(x_k,1)=a_{0}\,x_k^{n},
\]
which implies that the level curve \(H = a_{0}x_k^{n}\) contains a closed orbit \(\gamma_k\) of type “\emph{\textbf{bb}}”, connecting the points \((x_k,0)\) and \((x_k,1)\) on the torus (due to the identification of the opposite edges of the fundamental domain \([0,1]\times[0,1]\)).

The non-degeneracy of each cycle \(\gamma_k\) follows from the fact that \(H\) is homogeneous of degree \(n\), as expressed by the identity
\[
x\,\partial_x H(x,y) + y\,\partial_y H(x,y) = n\,H(x,y).
\]
Since along \(\gamma_k\) we have \(H(x_k,y) = a_0 x_k^n \neq 0\), it follows that \(\nabla H \neq 0\) throughout the orbit, ensuring a positive lower bound for \(\|\nabla H\|\).

Transversality of \(\gamma_k\) with respect to the “top” and “bottom” edges of the square is verified by evaluating the vector field at \((x_k,0)\) and \((x_k,1)\):
\[
X\cdot(0,1)\bigm|_{(x_k,0)}
=
-\partial_x H(x_k,0)
=
-n\,a_{0}\,x_k^{\,n-1}
\neq0,
\]
\[
X\cdot(0,-1)\bigm|_{(x_k,1)}
=
\partial_x H(x_k,1)
\neq0.
\]
These inequalities guarantee that the intersections with the boundary are transversal, i.e., the flow crosses the edges non-tangentially.

Finally, the sign choice in \(a_{0}\) ensures that each cycle \(\gamma_k\) remains entirely contained in \((0,1)\times[0,1]\), and the distinct level values \(a_{0}x_k^{n}\) ensure that these cycles are mutually isolated. Therefore, we conclude that there are exactly \(n-1\) limit cycles of type \emph{\textbf{bb}} located along the lines \(x = k/n\), which satisfy all the conditions of normalization, existence, non-degeneracy, transversality, and isolation required by the main theorem.

\end{exemplo}

\section*{Acknowledgements}

R. M. Martins was partially supported by FAPESP grants 2021/08031-9 and 2018/03338-6, CNPq grants 315925/2021-3, 434599/2018-2 and 306287/2024-2. This study was financed in part by the Coordenação de Aperfeiçoamento de Pessoal de Nível Superior - Brasil (CAPES) - Finance Code 001.


\bibliographystyle{unsrt}  
\bibliography{bibliografia.bib}

\begin{thebibliography}{10}

\bibitem{makarenkov2012dynamics}
Oleg Makarenkov and Jeroen~SW Lamb.
\newblock Dynamics and bifurcations of nonsmooth systems: A survey.
\newblock {\em Physica D: Nonlinear Phenomena}, 241(22):1826--1844, 2012.

\bibitem{teixeira2022perturbation}
Marco~Antonio Teixeira.
\newblock Perturbation theory for non-smooth systems.
\newblock In {\em Perturbation Theory: Mathematics, Methods and Applications}, pages 503--517. Springer, 2022.

\bibitem{bernardo2008piecewise}
Mario Bernardo, Chris Budd, Alan~Richard Champneys, and Piotr Kowalczyk.
\newblock {\em Piecewise-smooth dynamical systems: theory and applications}, volume 163.
\newblock Springer Science \& Business Media, 2008.

\bibitem{simpson2010bifurcations}
David John~Warwick Simpson.
\newblock {\em Bifurcations in piecewise-smooth continuous systems}, volume~70.
\newblock World Scientific, 2010.

\bibitem{guardia2011generic}
Marcel Guardia, TM~Seara, and Marco~Antonio Teixeira.
\newblock Generic bifurcations of low codimension of planar filippov systems.
\newblock {\em Journal of Differential equations}, 250(4):1967--2023, 2011.

\bibitem{broucke2001structural}
Mireille~E Broucke, C~Pugh, and Slobodan~N Simic.
\newblock Structural stability of piecewise smooth systems.
\newblock {\em Computational and applied mathematics}, 20(1-2):51--89, 2001.

\bibitem{jacquemard2011piecewise}
A~Jacquemard, MA~Teixeira, and DJ~Tonon.
\newblock Piecewise smooth reversible dynamical systems at a two-fold singularity.
\newblock {\em arXiv preprint arXiv:1111.4973}, 2011.

\bibitem{jacquemard2013stability}
ALAIN Jacquemard, MARCO~A Teixeira, and DURVAL~J Tonon.
\newblock Stability conditions in piecewise smooth dynamical systems at a two-fold singularity.
\newblock {\em Journal of Dynamical and Control Systems}, 19:47--67, 2013.

\bibitem{kuznetsov2003one}
Yu~A Kuznetsov, Sergio Rinaldi, and Alessandra Gragnani.
\newblock One-parameter bifurcations in planar filippov systems.
\newblock {\em International Journal of Bifurcation and chaos}, 13(08):2157--2188, 2003.

\bibitem{teixeira1990stability}
Marco~Antonio Teixeira.
\newblock Stability conditions for discontinuous vector fields.
\newblock {\em Journal of Differential Equations}, 88(1):15--29, 1990.

\bibitem{teixeira1993generic}
Marco~Antonio Teixeira.
\newblock Generic bifurcation of sliding vector fields.
\newblock {\em Journal of mathematical analysis and applications}, 176(2):436--457, 1993.

\bibitem{di2011two}
Mario di~Bernardo, Alessandro Colombo, and Enric Fossas.
\newblock Two-fold singularity in nonsmooth electrical systems.
\newblock In {\em 2011 IEEE International Symposium of Circuits and Systems (ISCAS)}, pages 2713--2716. IEEE, 2011.

\bibitem{chillingworth2002discontinuity}
DRJ Chillingworth.
\newblock Discontinuity geometry for an impact oscillator.
\newblock {\em Dynamical Systems}, 17(4):389--420, 2002.

\bibitem{jacquemard2012coupled}
Alain Jacquemard and Durval~J Tonon.
\newblock Coupled systems of non-smooth differential equations.
\newblock {\em Bulletin des Sciences Math{\'e}matiques}, 136(3):239--255, 2012.

\bibitem{jeffrey2009two}
Mike~R Jeffrey and Alessandro Colombo.
\newblock The two-fold singularity of discontinuous vector fields.
\newblock {\em SIAM Journal on Applied Dynamical Systems}, 8(2):624--640, 2009.

\bibitem{filippov1960differential}
Aleksei~Fedorovich Filippov.
\newblock Differential equations with discontinuous right-hand side.
\newblock {\em Matematicheskii sbornik}, 93(1):99--128, 1960.

\bibitem{llibre2018limit}
Jaume Llibre, Ricardo~Miranda Martins, and Durval~Jos{\'e} Tonon.
\newblock Limit cycles of piecewise smooth differential equations on two dimensional torus.
\newblock {\em Journal of dynamics and differential equations}, 30(3):1011--1027, 2018.

\bibitem{martins2019chaotic}
Ricardo~M Martins and Durval~J Tonon.
\newblock The chaotic behaviour of piecewise smooth differential equations on two-dimensional torus and sphere.
\newblock {\em Dynamical Systems}, 34(2):356--373, 2019.

\bibitem{griffiths1994principles}
Phillip Griffiths and Joseph Harris.
\newblock {\em Principles of Algebraic Geometry}.
\newblock Wiley-Interscience, New York, 1994.
\newblock Theorem I.3.1 on Bézout’s theorem.

\end{thebibliography}

\end{document}